\documentclass[11pt,twoside]{article}
\usepackage{amsmath, amsfonts, amssymb}
\usepackage[russian,english]{babel}
\usepackage{amscd,amssymb,amsmath}
\usepackage[arrow,matrix]{xy}
\usepackage{graphicx}
\usepackage{cite}
\usepackage{ifpdf}
\usepackage{tikz}
\begin{document}
\begin{center}\Large{\textbf{Classification of two-dimensional real evolution algebras and dynamics of some two-dimensional chains of evolution algebras}}

\textbf{Murodov Sh. N.}

e-mail: murodovs@yandex.ru

\textit{Institute of Mathematics. Tashkent}

\end{center}

{\bf Abstract.} In this paper we give classification of two-dimensional real evolution algebras.
For several chains of evolution algebras we study their classification dynamics.\\

{\it AMS classifications (2010):} 17D92; 37C99; 60J25.\\[2mm]

In [1] a notion of evolution algebra is introduced. This
evolution algebra is defined as follows.
 Let $(E,\cdot)$ be an algebra over a field $K$. If it admits a
basis $e_1,e_2,\dots$, such that $e_i\cdot e_j=0$, if $i\ne j$ and $
e_i\cdot e_i=\sum_{k}a_{ik}e_k$, for any $i$, then this algebra is
called an {\it evolution algebra}. This basis called natural basis of this algebra.

In this paper we consider finite dimensional evolution algebra $E$ over the field $\mathbb{R}$.

From the definition of the evolution algebra it is easy to see that this algebra is commutative, but
not associative, in general. It is important to note that, there exist several classes of non-associative algebras (baric,
evolution, Bernstein, train, stochastic, etc.), whose investigation
has provided a number of significant contributions to theoretical
population genetics. Such classes have been defined different times
by several authors, and all algebras belonging to these classes are
generally called "genetic"  [1, 4-9].

Following [2] we consider a family $\left\{E^{[s,t]}:\ s,t \in\mathbb{R},\ 0\leq s\leq t
\right\}$ of $n$-dimensional evolution algebras over the field $\mathbb{R}$,
with basis $e_1,\dots,e_n$ and multiplication table
\begin{equation}\label{1}
 e_ie_i =\sum_{j=1}^na_{ij}^{[s,t]}e_j, \ \ i=1,\dots,n; \ \
e_ie_j =0,\ \ i\ne j.\end{equation} Here parameters $s,t$ are
considered as time.

Denote by
$\mathcal{M}^{[s,t]}=\left(a_{ij}^{[s,t]}\right)_{i,j=1,\dots,n}$-the matrix
of structural constants.

{\bf Definition 1.} A family $\left\{E^{[s,t]}:\ s,t \in \mathbb{R},\ 0\leq s\leq t
\right\}$ of $n$-dimensional evolution algebras over the field $\mathbb{R}$
is called a chain of evolution algebras (CEA) if the matrix
$\mathcal{M}^{[s,t]}$ of structural constants satisfies the
Chapman-Kolmogorov equation
$\mathcal{M}^{[s,t]}=\mathcal{M}^{[s,\tau]}\mathcal{M}^{[\tau,t]}, \ \ \mbox{for any} \ \
s\leq \tau\leq t.$

Denote by
$\mathcal{T}=\left\{(s,t):0\leq s\leq t\right\}.$

{\bf Definition 2.} A CEA is called a time-homogenous CEA if
the matrix $\mathcal{M}^{[s,t]}$ depends only on $t-s$. In this case we write
$\mathcal{M}^{[t-s]}$.

{\bf Definition 3.} A CEA is called periodic if its matrix $\mathcal{M}^{[s,t]}$
is periodic with respect to at least one of the variables  $s$, $t$,
i.e. (periodicity with respect to $t$) $\mathcal{M}^{[s,t+P]}=\mathcal{M}^{[s,t]}$ for
all values of $t$. The constant $P$ is called the period, and is
required to be nonzero.

In [2, 10] some classes of chains of evolution algebras are studied.
Here we give matrices of structural constants of two dimensional chains of evolution algebras, which are constructed in [2]:
\begin{equation}\label{th}
\mathcal{M}_{1}^{[s,t]}={1\over 2}\left(\begin{array}{cc}
\lambda^t+\mu^t & \lambda^t-\mu^t\\[2mm]
\lambda^t-\mu^t & \lambda^t+\mu^t\\[2mm]
\end{array}\right),
\end{equation}
where $\lambda, \mu \geq 0$.

\begin{equation}\label{cs}
\mathcal{M}_2^{[s,t]}={1\over 2}\left(\begin{array}{cc}
\cos(t-s) & \sin(t-s)\\[2mm]
-\sin(t-s) & \cos(t-s)
\end{array}\right),
\end{equation}

\begin{equation}\label{ph}
\mathcal{M}_3^{[s,t]}={1\over 2}\begin{array}{ll}
\left(\begin{array}{cc}
1+\Phi(t)(\Psi(t)-\Psi(s))+{\Phi(t)\over \Phi(s)} & 1-\Phi(t)(\Psi(t)-\Psi(s))-{\Phi(t)\over \Phi(s)}\\
1+\Phi(t)(\Psi(t)-\Psi(s))-{\Phi(t)\over \Phi(s)} & 1-\Phi(t)(\Psi(t)-\Psi(s))+{\Phi(t)\over \Phi(s)}
\end{array}\right)
\end{array},
\end{equation}
where $\Phi(s)\neq 0, \Psi(s)$ are arbitrary functions.

We denote by $E_i^{[s,t]}, 0\leq s\leq t$ chains of evolution algebras corresponding
to $\mathcal{M}_i^{[s,t]}, i=1,2,3.$

We note that the algebra $E_1^{[s,t]}$ corresponding to (\ref{th}) is time-homogenous CEA,
the algebra $E_2^{[s,t]}$ corresponding to (\ref{cs}) is
a periodic CEA and the algebra $E_3^{[s,t]}$ is time non-homogenous CEA.

In [3] a classification of two dimensional complex evolution algebras is given.

To study dynamics of chains of evolution algebras we need the next theorem,
which gives the classification of two-dimensional real evolution algebras.

\textbf{Theorem 1.} Any two-dimensional real evolution algebra \textit{E} is isomorphic
to one of the following pairwise non-isomorphic algebras:

(i) $dimE^2=1$

$E_{1}: e_1e_1=e_1, \ \ e_2e_2=0;$

$E_{2}: e_1e_1=e_1, \ \  e_2e_2=e_1;$

$E_{3}: e_1e_1=e_1+e_2, \ \  e_2e_2=-e_1-e_2;$

$E_4: e_1e_1=e_2, \ \  e_2e_2=0 ;$

$E_5: e_1e_1=e_2, \ \  e_2e_2=-e_2;$

(ii) $dimE^2=2$:

$E_6(a_2;a_3): e_1e_1=e_1+a_2 e_2, \ \ e_2e_{2}=a_3 e_1+e_2, \ \ 1-a_2 a_3 \neq 0, a_2, a_3\in {\mathbb{R}}.$ Moreover $E_6(a_2;a_3)$
 is isomorphic to $E_6(a_3;a_2).$

$E_7(a_4): e_1e_1=e_2, \ \ e_2e_{2}=e_1+a_4 e_2, \ \ where \ a_4\in {\mathbb{R}}.$

\textit{Proof}. For a general two-dimensional evolution algebra we have $e_1e_1=a_1e_1+a_2e_2, \ \
 e_2e_2=a_3e_1+a_4e_2$ and $e_1e_2=e_2e_1=0, a_i\in\mathbb{R}.$

\ \ \textbf{(i)} Since $dimE^2=1,$ we have $e_1e_1=c_1(a_1e_1+a_2e_2), \ \ e_2e_2=c_2(a_1e_1+a_2e_2)$ and $e_1e_2=e_2e_1=0.$
Evidently, $(c_1,c_2)\neq (0,0)$, because otherwise our algebra will be abelian. Since $e_1$ and $e_2$ are
symmetric, we can suppose $c_1 \neq0,$ and by a simple change of basis we can suppose $c_1=1.$

\textbf{Case 1.} $a_1\neq0.$ We take an appropriate change of basis $e'_1=a_1e_1+a_2e_2, \ e'_2=Ae_1+Be_2,$ where $a_1B-a_2A\neq0.$
Consider the product

$0=e'_1 e'_2=(a_1e_1+a_2e_2)(Ae_1+Be_2)=a_1A(a_1e_1+a_2e_2)+a_2Bc_2(a_1e_1+a_2e_2)=\\
=(a_1A+a_2Bc_2)(a_1e_1+a_2e_2)$

Therefore, $a_1A+a_2Bc_2=0,$ i.e., $A=-\frac{a_2Bc_2}{a_1}$ and $a_1B-a_2A=a_1B+\frac{a^2_2Bc_2}{a_1}\neq0.$

For $B\neq 0,$ it means that in the case when $a^2_1+a^2_2c_2\neq0$ we can take the above change.

Consider the products

$e'_1e'_1=(a_1e_1+a_2e_2)(a_1e_1+a_2e_2)=a^2_1(a_1e_1+a_2e_2)+a^2_2c_2(a_1e_1+a_2e_2)=\\
=(a^2_1+a^2_2c_2)(a_1e_1+a_2e_2)=(a^2_1+a^2_2c_2)e'_1,$

$e'_2e'_2=(Ae_1+Be_2)(a_1e_1+a_2e_2)=A^2(a_1e_1+a_2e_2)+B^2c2(a_1e_1+a_2e_2)=\\
=(A^2+B^2c_2)(a_1e_1+a_2e_2)=(\frac{a^2_2B^2c^2_2}{a^2_1}+B^2c_2)e'_1=\frac{B^2c_2(a^2_1+a^2_2c_2)}{a^2_1}e'_1.$

 \  \ \textbf{Case 1.1.}  $c_2=0.$ Then $e_1e_1=a^2_1e_1$ and $e_2e_2=e_1e_2=e_2e_1=0.$ Taking

$ e'_1=\frac{e_1}{a^2_1}$, we obtain the algebra $E_1.$

 \ \\textbf{ Case 1.2.} $c_2\neq0$. Then taking $B=\sqrt{\frac{a^2_1}{|c_2|}}$, we obtain $e_1e_1=(a^2_1+a^2_2c_2)e_1, \ \ e_2e_2=(a^2_1+a^2_2c_2)e_1,$
 when $c_2>0$ and also $e_1e_1=(a^2_1+a^2_2c_2)e_1, \ \ e_2e_2=-(a^2_1+a^2_2c_2)e_1,$ when $c_2<0$.

 If $c_2>0$ then $a^2_1+a^2_2c_2\neq0$ and we can take the change of basis $e'_1=\frac{e_1}{a^2_1+a^2_2c_2}, \ \ e'_2=\frac{e_2}{a^2_1+a^2_2c_2}$
 which gives the algebra $E_2$ with multiplication table $e_1e_1=e_1, \ \ e_2e_2=e_1.$

If $a^2_1+a^2_2c_2\neq0$ when $c_2<0$ and we can take the change of basis $e'_1=\frac{e_1}{a^2_1+a^2_2c_2}, \ \ e'_2=-\frac{e_2}{a^2_1+a^2_2c_2}$
 which gives the algebra with multiplication table $e_1e_1=e_1, \ \ e_2e_2=-e_1.$ It is easy to check that this algebra is isomorphic
 to the algebra $E_5$ with the change of basis $e'_1=-e_2, \ e'_2=e_1.$

 If $a^2_1+a^2_2c_2=0$ $(c_2<0)$, since $a_1\neq0$ we have $a_2\neq0$ then $c_2=-\frac{a^2_1}{a^2_2}$ and we have $e_1e_1=a_1e_1+a_2e_2$
 and $e_2e_2=-\frac{a^3_1}{a^2_2}e_1-\frac{a^2_1}{a_2}e_2.$ Then the change of basis $e'_1=\frac{e_1}{a_1}, e_2=\frac{a_2}{a^2_1}e_2$ gives
  the algebra $E_3.$

  \textbf{Case 2. }$a_1=0.$ Then we have $e_1e_1=a_2e_2$ and $e_2e_2=c_2a_2e_2,$ where $a_2\neq0$.

  If $c_2=0$ then by the change $e'_1=\frac{e_1}{\sqrt{|a_2|}}$ we get the algebra $E_4$ when $a_2>0.$
  When $a_2<0$ by this change of basis we get the algebra $e_1e_1=-e_2, e_2e_2=0$ which is
  isomorphic to the algebra $E_4$ by the change of basis $e'_1=e_1, e'_2=e_2$.

  If $c_2\neq0$, then by $e'_1=\frac{e_1}{\sqrt{|c_2|a^2_2}}$ and $e'_2=\frac{e_2}{c_2a_2},$
  we get the algebra $e_1e_1=e_2, \ \ e_2e_2=e_2 \ (c_2>0)$ which is isomorphic to the algebra $E_2$.

  If $c_2<0$ then by the change of basis $e'_1=\frac{e_1}{\sqrt{|c_2|a^2_2}}$ and $e'_2=\frac{e_2}{c_2a_2},$
  we will take the algebra with multiplication table  $e_1e_1=-e_2, \ \ e_2e_2=e_2$ which is isomorphic to $E_5$.

  \textbf{(ii)} Now consider algebras with $dimE^2=2.$ Let us write $e_1e_1=a_1e_1+a_2e_2, \ \ e_2e_2=a_3e_1+a_4e_2,$
  where $a_1a_4-a_2a_3\neq0.$

  \textbf{Case 1.} $a_1\neq0$ and $a_4\neq0$. Then the change of basis $e_1=a^{-1}_1e_1, e_2=a^{-1}_4e_2$ makes
  possible to suppose $a_1=a_4=1$. Therefore, we have two-parametric family $E_7(a_2,a_3): e_1e_1=e_1+a_2e_2, \ \
  e_2e_2=a_3e_1+e_2, \ \ 1-a_2a_3\neq0.$
  Let us take the general change of basis $e'_1=A_1e_1+A_2e_2, \ \ e'_2=B_1e_1+B_2e_2,$ where $A_1B_2-A_2B_1\neq0.$
  Consider the product

  $0=e'_1e'_2=(A_1e_1+A_2e_2)(B_1e_1+B_2e_2)=A_1B_1(e_1+a_2e_2)+A_2B_2(a_3e_1+e_2)
  =(A_1B_1+A_2B_2a_3)e_1+(A_1B_1a_2+A_2B_2)e_2.$

  Since in this new basis the algebra should be also an evolution algebra, we have
  $A_1B_1+A_2B_2a_3=0$ and $A_1B_1a_2+A_2B_2=0.$ From this we have $A_2B_2(1-a_2a_3)=0$
  and $A_1B_1(1-a_2a_3)=0.$ Since $1-a_2a_3\neq0$, we have $A_1B_1=A_2B_2=0.$

  \textbf{Case 1.1.} Let $A_2=0$. Then $B_1=0$. Consider the products

  $e'_1e'_1=A^2_1(e_1+a_2e_2)=e'_1+a'_2e'_2=A_1e_1+a'_2B_2e_2$

            $\Rightarrow A^2_1=A_1, A^2_1a_2=a'_2B_2 \Rightarrow A_1=1,$

  $e'_2e'_2=B^2_2(a_3e_1+e_2)=a'_3e'_1+e'_2=a'_3A_1e_1+B_2e_2$

            $\Rightarrow B^2_2a_3=a'_3A_1, B^2_2=B_2 \Rightarrow B_2=1.$

  \textbf{Case 1.2.} Let $A_1=0$. Then $B_2=0$, and from the family of algebras $E_6(a_2,a_3)$ we get the family $E_6(a_3,a_2)$.

  Case 2. Let $a_1=0$ or $a_4=0$. Since $e_1$ and $e_2$ are symmetric, without loss of generality we can suppose
  $a_1=0$, i.e., $e_1e_1=a_2e_2$ and $e_2e_2=a_3e_1+a_4e_2,$ where $a_2a_3\neq0.$
  Taking the change of basis $e'_1=\sqrt[3]{\frac{1}{a^2_2a_3}}e_1, \ \ e'_2=\sqrt[3]{\frac{1}{a_2a^2_3}}e_2$, we obtain
  one of the parametric family of algebras $E_7(a_4): \ e_1e_1=e_2, \ \ e_2e_2=e_1+a_4e_2.$ Theorem is proved.\\

\textbf{Remark.} We note that the classification of two dimensional complex evolution algebras consists complex variant of
algebras $E_i, \ i=1,2,3,5,6,7$ [3]. But $E_4$ is present only in real case.

For studying the dynamics of the CEA listed above we need the next lemma.\\

\textbf{Lemma 1.} Evolution algebra corresponding to a matrix

$$\textrm{(i)} \ \ \ \left(\begin{array}{cc}
\lambda & \mu\\
\mu & \lambda
\end{array}\right) \, \simeq \, \left\{\begin{array}{ll}
\begin{array}{cc}
E_0, \ \ \mbox {if} \ \ \lambda=\mu=0\end{array};\\[4mm]
\begin{array}{cc}
E_2, \ \ \mbox {if} \ \ \lambda=\mu\neq0\end{array};\\[4mm]
\begin{array}{ll}
E_6(\frac{\mu}{\lambda};\frac{\mu}{\lambda}), \ \ \mbox{if} \ \ \lambda\neq\mu, \lambda\neq0, \mu\in\mathbb{R}
\end{array};\\[4mm]
\begin{array}{ll}
E_7(0), \ \ \mbox{if} \ \lambda\neq\mu, \lambda=0, \mu\neq 0
\end{array};
\end{array}\right.$$

$$\textrm{(ii)} \ \ \ \left(\begin{array}{cc}
\lambda & \mu\\
-\mu & \lambda
\end{array}\right) \ \ \mbox{is isomorphic to} \left\{\begin{array}{ll}
\begin{array}{cc}
E_0, \ \ \mbox{if} \ \ \lambda=\mu=0
\end{array};\\[4mm]
\begin{array}{cc}
E_6(\frac{\mu}{\lambda};-\frac{\mu}{\lambda}), \ \ \mbox{if} \ \ \lambda\neq0, \mu\in{\mathbb{R}}
\end{array};\\[4mm]
\begin{array}{cc}
E_7(0), \ \ \mbox{if} \ \ \lambda=0, \mu\neq0
\end{array};\\[4mm]
\end{array}\right.$$

$$\textrm{(iii)} \ \ \ \left(\begin{array}{cc}
1+\lambda & 1-\lambda\\
1+\mu & 1-\mu
\end{array}\right) \ \ \mbox{is isomorphic to} \left\{\begin{array}{ll}
\begin{array}{cc}
E_2, \ \ \mbox{if} \ \ \lambda=\mu
\end{array};\\[4mm]
\begin{array}{cc}
E_6(\frac{(1+\lambda)(1+\mu)}{(1-\mu)^2};\frac{(1-\lambda)(1-\mu)}{(1+\lambda)^2}), \ \
 \mbox{if} \ \ \lambda\neq\mu \ \mbox{and} \ \lambda\neq-1, \mu\neq1
\end{array};\\[4mm]
\begin{array}{cc}
E_7(\frac{1-\mu}{\sqrt[3]{2(1+\mu^2)}}), \ \ \mbox{if} \ \ \lambda=-1, \mu\neq-1 \end{array};\\[4mm]
\begin{array}{cc}
E_7(\frac{\lambda}{\sqrt[3]{(1-\lambda)^2}}), \ \ \mbox{if} \ \ \lambda\neq1, \mu=1 \end{array},
\end{array}\right.$$
where $E_0$ is the trivial evolution algebra (i.e., with zero multiplication) and
evolution algebras $E_i, \ i=\overline{1,7}$ are given in Theorem 1.

\textbf{Proof}. Let
$$\mathcal{M}=\left(\begin{array}{cc}
\alpha  & \beta \\[2mm]
\gamma  & \delta
\end{array}\right), \ \ \mathcal{A}=\left(\begin{array}{cc}
a & b\\[2mm]
c & d
\end{array}\right)$$
be matrices of structural constants of evolution algebras $E_{\mathcal{M}}$ and
$E_{\mathcal{A}}$. The multiplication table in $E_{\mathcal{M}}$ is

$$e_1e_1=\alpha e_1+\beta e_2, \ \ e_2e_2=\gamma e_1+\delta e_2$$

and in $E_{\mathcal{A}}$ is

$$e'_1e'_1=ae'_1+be'_2, \ \  e'_2e'_2=ce'_1+de'_2.$$

Let

$$e'_1=xe_1+ye_2, \ \  e'_2=ze_1+ve_2$$

be change of basis, where $xv-yz\neq 0.$

We have the following equations: \\
$0=e'_1e'_2=(xe_1+ye_2)(ze_1+ve_2)=(\alpha xz+\gamma yv)e_1
+(\beta xz+\delta yv)e_2$,

$e'_1e'_1=(xe_1+ye_2)(xe_1+ye_2)=(\alpha x^2+\gamma y^2)e_1+
(\beta x^2+\delta y^2)e_2,$

$e'_1e'_1=a(xe_1+ye_2)+b(ze_1+ve_2)=
(ax+bz)e_1+(ay+bv)e_2.$

$e'_2e'_2=(ze_1+ve_2)(ze_1+ve_2)=(\alpha z^2+\gamma v^2)e_1+
(\beta z^2+\delta v^2)e_2,$

$e'_2e'_2=c(xe_1+ye_2)+d(ze_1+ve_2)=
(cx+dz)e_1+(cy+dv)e_2.$

Consequently,
\begin{equation}\label{s1}
\left\{\begin{array}{lllllll}
xv-yz\neq0 \\
\alpha xz+\gamma yv=0 \\
\beta xz+\delta yv=0 \\
\alpha x^2+\gamma y^2=ax+bz \\
\beta x^2+\delta y^2=ay+bv \\
\alpha z^2+\gamma v^2=cx+dz \\
\beta z^2+\delta v^2=cy+dv
\end{array}\right..
\end{equation}
Therefore, we should solve system of equations (\ref{s1}) for each evolution algebra.
Since CEAs listed above have matrices as in Lemma 1, it will be enough to check only these forms of matrices.

\textbf{Case \textrm{(i)}.} Let $\alpha =\delta =\lambda , \beta =\gamma=\mu $.
Then (\ref{s1}) will be

$$\left\{\begin{array}{lllllll}
xv-yz\neq0 \\
\lambda xz+\mu yv=0 \\
\mu xz+\lambda yv=0 \\
\lambda x^2+\mu y^2=ax+bz \\
\mu x^2+\lambda y^2=ay+bv \\
\lambda z^2+\mu v^2=cx+dz \\
\mu z^2+\lambda v^2=cy+dv
\end{array}\right..$$

\textbf{Case \textrm{(i)}.1.} For $\lambda=\mu=0$ this algebra will be trivial EA.

\textbf{Case \textrm{(i)}.2.} For $\lambda=\mu\neq0$ this algebra will be isomorphic to $E_2$ by the change of basis
$e'_1=\frac {1}{2\lambda }e_1+\frac {1}{2\lambda }e_2, e'_2=-\frac {1}{2\lambda }e_1+\frac {1}{2\lambda }e_2$.

\textbf{Case \textrm{(i)}.3.} Let $\lambda=0, \lambda\neq\mu, \mu \neq 0$ then this algebra will be isomorphic to $E_7(0)$ by the change of basis $e'_1=\frac {1}{\mu}e_1, e'_2=\frac{1}{\mu}e_2.$

\textbf{Case \textrm{(i)}.4.} Let $\lambda\neq 0$, then this algebra will be isomorphic to $E_6(\frac{\mu}{\lambda},\frac{\mu}{\lambda})$
by change of basis $e'_1=\frac {1}{\lambda}e_1, e'_2=\frac{1}{\lambda}e_2.$

\textbf{Case \textrm{(ii)}.} Let $\alpha=\delta=\lambda, \ \beta=-\gamma=\mu$.
Then (\ref{s1}) will be

$$\left\{\begin{array}{lllllll}
xv-yz\neq0 \\
\lambda xz-\mu yv=0 \\
\mu xz+\lambda yv=0 \\
\lambda x^2-\mu y^2=ax+bz \\
\mu x^2+\lambda y^2=ay+bv \\
\lambda z^2-\mu v^2=cx+dz \\
\mu z^2+\lambda v^2=cy+dv
\end{array}\right..$$

\textbf{Case \textrm{(ii)}.1.} For $\lambda=0, \mu\neq 0$ this algebra will be isomorphic to
$E_7(0)$ by the change of basis $e'_1=-\frac{1}{\mu}e_1, e'_2=\frac{1}{\mu}e_2.$

\textbf{Case \textrm{(ii)}.2.} For $\lambda\neq 0 \ \ \mbox{and} \ \ \mu \in\mathbb{R}$ this algebra will be isomorphic to
$E_6(\frac{\mu}{\lambda};-\frac{\mu}{\lambda})$
by the change of basis $e'_1=\frac{1}{\lambda}e_1, e'_2=\frac{1}{\lambda}e_2.$

\textbf{Case \textrm{(iii)}.} Let $\alpha=1+\lambda, \beta=1-\lambda, \gamma=1+\mu, \delta=1-\mu$.
Then (\ref{s1}) will be

\begin{equation}
\left\{\begin{array}{lllllll}
xv-yz\neq0 \\
(1+\lambda)xz+(1+\mu)yv=0 \\
(1-\lambda)xz+(1-\mu)yv=0 \\
(1+\lambda)x^2+(1+\mu)y^2=ax+bz \\
(1-\lambda)x^2+(1-\mu)y^2=ay+bv \\
(1+\lambda)z^2+(1+\mu)v^2=cx+dz \\
(1-\lambda)z^2+(1-\mu)v^2=cy+dv
\end{array}\right..
\end{equation}

\textbf{Case \textrm{(iii)}.1.} Case $\lambda=\mu$ this algebra will be isomorphic to $E_2$ by the change of basis
$e'_1=\frac{1+\lambda}{2(1+\lambda^2)}e_1+\frac{1-\lambda}{2(1+\lambda^2)}e_2,
e'_2=\frac{1-\lambda}{2(1+\lambda^2)}e_1-\frac{1+\lambda}{2\lambda^2-2\lambda+1}e_2.$

\textbf{Case \textrm{(iii)}.2.} For $\lambda\neq-1, \mu\neq1$ this algebra
will be isomorphic to $E_6(\frac{(1+\lambda)(1+\mu)}{(1-\mu)^2};\frac{(1-\lambda)(1-\mu)}{(1+\lambda)^2})$
by the change of basis $e'_1=\frac{1}{1-\mu}e_2, e'_2=\frac{1}{1+\lambda}e_1.$

\textbf{Case \textrm{(iii)}.3.} For $\lambda=-1, \mu\neq-1$ this algebra will be isomorphic
to $E_7(\frac{1-\mu}{\sqrt[3]{2(1+\mu^2)}})$
by the change of basis $e'_1=\frac{1}{\sqrt[3]{4(1+\mu)}}e_1, e'_2=\frac{1}{\sqrt[3]{2(1+\mu)^2}}e_2.$

\textbf{Case \textrm{(iii)}.5.} For $\lambda\neq1, \mu=1$ this algebra will be isomorphic
to $E_7(\frac{\lambda}{\sqrt[3]{(1-\lambda)^2}})$
by the change of basis $e'_1=\frac{1}{\sqrt[3]{4(1-\lambda)}}e_2, e'_2=\frac{1}{\sqrt[3]{2(1-\lambda)^2}}e_1.$

The next theorem gives the dynamics of chains of evolution algebras listed above, which we
shall prove by the above lemma.

{\bf Theorem 2.}

$$ E_{1}^{[s,t]} \ \ \mbox{in any time is isomorphic to} \ \
E_6\left(\frac{\lambda^t-\mu^t}{\lambda^t+\mu^t};\frac{\lambda^t-\mu^t}{\lambda^t+\mu^t}\right). \ \
\mbox{for all} \ \ \lambda, \mu;$$

$$E_2^{[s,t]} \simeq \left\{\begin{array}{ll}
\begin{array}{ll}
E_6(\tan(t-s);-\tan(t-s))\\ \ \ \mbox{for all} \ \ (s,t)\in \left\{(s,t):t\neq s+\frac{\pi}{2}+\pi k, \ \ k\in \mathbb{Z}\right\}
\end{array};\\[7mm]
\begin{array}{cc}
E_7(0) \ \ \mbox{for all} \ \ (s,t)\in \left\{(s,t): t=s+\frac{\pi}{2}+\pi k, \ \ k\in \mathbb{Z}\right\}\end{array};
\end{array}\right.$$

$$E_{3}^{[s,t]} \simeq \left\{\begin{array}{ll}
\begin{array}{l}
E_6\left(\frac{(1+\xi)(1+\zeta)}{(1-\zeta)^2};\frac{(1-\xi)(1-\zeta)}{(1+\xi)^2}\right)
 \ \ \mbox{for all} \ \ \xi\neq\zeta \ \mbox{and} \ \xi\neq-1, \zeta\neq1
\end{array};\\[7mm]
\begin{array}{l}
E_7\left(\frac{1-\zeta}{\sqrt[3]{2(1+\zeta^2)}}\right)
 \ \ \mbox{if}  \ \ \xi=-1, \zeta\neq-1
\end{array};\\[7mm]
\begin{array}{l}
E_7\left(\frac{\xi}{\sqrt[3]{(1-\xi)^2}}\right)
 \ \ \mbox{if} \ \ \xi\neq1, \zeta=1
\end{array},\\[7mm]
\end{array}\right.$$
where $\xi=1+\Phi(t)(\Psi(t)-\Psi(s))+{\Phi(t)\over \Phi(s)}, \zeta=1+\Phi(t)(\Psi(t)-\Psi(s))-{\Phi(t)\over \Phi(s)}$.

\textbf{Proof}. From the proved lemma it is easy to see that $E_1^{[s,t]}$ is isomorphic to
$E_6(\frac{\lambda^t-\mu^t}{\lambda^t+\mu^t};\frac{\lambda^t-\mu^t}{\lambda^t+\mu^t})$, by change of basis
$e'_1=\frac{1}{2}\frac{1}{\lambda^t+\mu^t}e_1, \ \ e'_2=\frac{1}{2}\frac{1}{\lambda^t+\mu^t}e_2.$

The CEA $E_2^{[s,t]}$ is isomorphic to $E_6(\tan(t-s);-\tan(t-s))$,
$\forall s,t\in\mathcal T,$ when $t\neq s+\frac{\pi}{2}+\pi k, k\in \mathbb{Z}$
by change of basis $e'_1=\frac{1}{\cos(t-s)}e_1, \ \ e'_2=\frac{1}{\cos(t-s)}e_2,$ and it is
isomorphic to $E_7(0)$ $\forall s,t\in\mathcal T, \ \ \mbox{when} \ \ t=s+\frac{\pi}{2}+\pi k, k\in \mathbb{Z}$
by change of basis $e'_1=-e_1, \ \ e'_2=e_2.$

The CEA $E_3^{[s,t]}$ is isomorphic to $E_6(\frac{(1+\xi)(1+\zeta)}{(1-\zeta)^2};\frac{(1-\xi)(1-\zeta)}{(1+\xi)^2})$ by the
change of basis $e'_1=\frac{1}{1-\zeta}e_2, e'_2=\frac{1}{1+\xi}e_1$, when $\xi\neq\zeta \ \mbox{and} \ \xi\neq-1, \zeta\neq1$.

When $\xi=1, \zeta\neq-1$ this CEA is isomorphic to $E_7(\frac{1-\zeta}{\sqrt[3]{2(1+\zeta^2)}})$ by the
change of basis $e'_1=\frac{1}{\sqrt[3]{4(1+\zeta)}}e_1, e'_2=\frac{1}{\sqrt[3]{2(1+\zeta)^2}}e_2$ and it is
isomorphic to $E_7(\frac{\xi}{\sqrt[3]{(1-\xi)^2}})$ by the change of basis
$e'_1=\frac{1}{\sqrt[3]{4(1-\xi)}}e_2, e'_2=\frac{1}{\sqrt[3]{2(1-\xi)^2}}e_1$, when
$\xi\neq1, \zeta=1$, where $\xi=1+\Phi(t)(\Psi(t)-\Psi(s))+{\Phi(t)\over \Phi(s)}, \zeta=1+\Phi(t)(\Psi(t)-\Psi(s))-{\Phi(t)\over \Phi(s)}$.
Theorem is proved.

\begin{center}\textbf{References}\end{center}

1. J. P. Tian, \emph{Evolution algebras and their applications},
Lecture Notes in Mathematics, 1921, Springer-Verlag, Berlin, 2008.\\
2. J.M. Casas, M. Ladra, U.A. Rozikov, \emph{A chain of evolution algebras}. Linear Algebra Appl.  435(4), 852--870  (2011).\\
3. J.M. Casas, M. Ladra, B.A. Omirov, U.A. Rozikov, \emph{On evolution algebras.} arXiv:1004.1050. To appear in Algebra Colloquium.\\
4. I.M.H. Etherington, \emph{Genetic algebras}, Proc. Roy. Soc. Edinburgh. 59, 242--258 (1939).\\
5. I.M.H. Etherington, \emph{Duplication of linear algebras}, Proc. Edinburgh Math. Soc. (2) 6, 222--230 (1941).\\
6. I.M.H. Etherington, \emph{Non-associative algebra and the simbolism of genetics},
Proc. Roy. Soc. Edinburgh. 61, 24--42 (1941).\\
7. Y.I. Lyubich, \emph{Mathematical structures in population
genetics}, Springer-Verlag, Berlin, 1992.\\
8. M.L. Reed, \emph{Algebraic structure of genetic inheritance},  Bull. Amer. Math.
Soc. (N.S.)  34(2), 107--130  (1997).\\
9. A. W\"orz-Busekros, \emph{Algebras in genetics}, Lecture Notes in
Biomathematics, 36. Springer-Verlag, Berlin-New York, 1980.\\
10. U.A. Rozikov, Sh.N. Murodov, Dynamics of two-dimensional evolution algebras. arXiv:1202.2690.\\

\end{document}